\documentclass[12pt]{amsart}
\usepackage{amsmath, amsfonts,amsthm,times,graphics}
\input amssym.def
\input amssym

\begin{document}

\newtheorem{theorem}{Theorem}[section]
\newtheorem{lemma}[theorem]{Lemma}
\newtheorem{corollary}[theorem]{Corollary}
\newtheorem{proposition}[theorem]{Proposition}
\newtheorem{conjecture}[theorem]{Conjecture}
\newtheorem{question}[theorem]{Question}
    \theoremstyle{definition}
\newtheorem{definition}[theorem]{Definition}
\newtheorem{example}[theorem]{Example}
\newtheorem{xca}[theorem]{Exercise}
\newtheorem{remark}[theorem]{{\it Remark}}

\newcommand{\abs}[1]{\lvert#1\rvert}

 \makeatletter

\title[On some metric topologies on Privalov spaces \ldots]
{On some metric topologies\\ on Privalov spaces on the unit disk}

\author{Romeo Me\v{s}trovi\'{c}${}^{\,1)}$}
\address{${}^{\,1)}$Maritime Faculty Kotor, University of Montenegro, 
85330 Kotor, Montenegro} \email{romeo@ucg.ac.me}

\author{\v{Z}arko Pavi\'{c}evi\'{c}${}^{\, 2),\, 3)}$}

\address{${}^{2)}$ Faculty of Science, University of Montenegro, 
81000  Podgorica, Montenegro; $\quad\quad\,\,$ ${}^{3)}$National Research Nuclear University MEPhI, 
Moscow Engineering Physics Institute, Moscow, Russia}
 \email{zarkop@ucg.ac.me}

\makeatother

{\renewcommand{\thefootnote}{}
\setcounter{footnote}{0}}
 \maketitle

\vspace*{7mm}

\noindent{\bf Summary.} 
Let $N^p$   $(1<p<\infty)$ be the  Privalov class $N^p$  of holomorphic 
functions  on the open unit disk $\Bbb  D$
in the complex plane. 
In 1977 M. Stoll proved that  the class $N^p$
equipped with the topology given by the metric $\lambda_p$ defined by 
 $$
\lambda_p(f,g)=\Bigg(\int_0^{2\pi}\big(\log(1+
\vert f^*(e^{i\theta})-g^*(e^{i\theta})\vert)\big)^p\,\frac{d\theta}
{2\pi}\Bigg)^{1/p},\quad f,g\in N^p,
 $$ 
becomes an $F$-algebra. In the recent overview paper by Me\v{s}trovi\'{c} and
  Pavi\'{c}evi\'{c} (2017) 
a survey of some known results on the 
topological structures of the Privalov spaces $N^p$ $(1<p<\infty)$
and their Fr\'{e}chet envelopes $F^p$ are presented.

In this  article we continue a survey of results concerning 
the topological structures of the spaces $N^p$ $(1(p<\infty)$.
In particular, for each $p>1$, we consider the  class
$N^p$ as the space $M^p$ equipped with the topology 
induced by the metric $\rho_p$ defined as
   $$
\rho_p(f,g)=\Bigg(\int_0^{2\pi}\log^p(1+M(f-g)(\theta))\,
\frac{d\theta}{2\pi}\Bigg)^{1/p},\quad f,g\in M^p,
  $$
where $Mf(\theta)=\sup_{0\leqslant r<1}\big\vert 
f\big(re^{i\theta})\big\vert.$

On the other hand, we consider 
the class $N^p$ with the metric topology 
introduced by Me\v{s}trovi\'{c},  Pavi\'{c}evi\'{c} 
and  Labudovi\'{c} (1999) which generalizes the 
Gamelin-Lumer's  metric which is generally defined on a measure space
$(\Omega, \Sigma, \mu)$ with a positive finite measure
$\mu$. The space $N^p$ with the associated modular in the sense 
of Musielak and Orlicz becomes the 
Hardy-Orlicz class.  
It is  noticed that the all considered metrics induce the same
  topology on the space $N^p$.

\vspace*{10mm}

\noindent {\bf 1 INTRODUCTION AND PRELIMINARY RESULTS}

\vspace{3mm}

Let $\Bbb D$ denote the open unit disk 
in the complex plane  and let $\Bbb T$ denote the
boundary of $\Bbb D$.  Let $L^q(\Bbb T)$ $(0<q\le \infty)$ be the
familiar Lebesgue spaces on the unit circle $\Bbb T$.
The {\it  Privalov 
class} $N^p$ $(1<p<\infty)$ consists of all holomorphic functions $f$ on $\Bbb D$ 
for  which 
   $$
\sup_{0\le r<1}\int_0^{2\pi}(\log^+\vert f(re^{i\theta})\vert)^p\,
\frac{d\theta}{2\pi}<+\infty, \eqno(1)
   $$ 
where for $z\in\Bbb C$, $\log^+|z|=\max(\log|z|,0)$ if  $z\not= 0$
 and $\log^+0=0$. These classes were firstly considered  by I.I. Privalov in 
\cite[p. 93]{p}, where $N^p$ is denoted as $A_q$. 

Notice that for $p=1$ the condition (1) defines the 
{\it Nevanlinna class} $N$ of holomorphic functions on $\Bbb D$.
 Recall that the {\it Smirnov class} $N^+$  is the set
of all functions  $f\in N$  such that
   $$
\lim_{r\rightarrow 1}\int_0^{2\pi}\log^+\vert f(re^{i\theta})\vert
\,\frac{d\theta}{2\pi}=\int_0^{2\pi}\log^+\vert f^*(e^{i\theta})\vert\,
\frac{d\theta}{2\pi}<+\infty,
   $$
where $f^*$ is the boundary function of $f$ on $\Bbb T$, i.e.,
 $$
f^*(e^{i\theta})=\lim_{r\rightarrow 1-}f(re^{i\theta})
 $$ is the {\it radial
 limit} of $f$  which exists for almost every $e^{i\theta}\in\Bbb T$. 
Recall that the classical {\it Hardy space} $H^q$ $(0<q\le\infty)$ consists of all
functions $f$ holomorphic on $\Bbb D$ such that 
\begin{equation*} 
\sup_{0\le r<1}\int_0^{2\pi}\big\vert f\big(re^{i\theta}\big)\big\vert^q
\frac{d\theta}{2\pi}<\infty
\end{equation*}
if $0<q<\infty$, and which are bounded when $q=\infty$:
\begin{equation*}
\sup_{z\in \Bbb D}\vert f(z)\vert<\infty{.}
\end{equation*}
It is known that (see  \cite{Mo} and \cite{mp8}) 
   $$
N^r\subset N^p\;(r>p),\quad\bigcup_{q>0}H^q\subset
\bigcap_{p>1}N^p,\quad{\rm and}\quad\bigcup_{p>1}N^p\subset M\subset N^+
\subset N,
   $$
where the above containment relations are proper.

It is well known (see, e.g., \cite[p. 26]{D}) that a function $f\in N^+$
has a unique factorization of the form
  $$
f(z)=B(z)S_{\mu}(z)F(z), \quad z\in \Bbb D, 
   $$
where $B$ is the {\it Blaschke product} with respect to zeros
$\{z_k\}\subset \Bbb D$ of $f$, $S_{\mu}$ is a {\it singular inner function} 
and $F$ is an {\it outer function} in $N^+$, i.e.,
 $$
B(z)=z^m\prod_{k=1}^\infty\frac{|z_k|}{z_k}\cdot
\frac{z_k-z}{1-\bar{z}_kz}, \quad z\in \Bbb D, 
 $$
with $\sum_{k=1}^\infty\big(1-|z_k|\big)<\infty$, $m$ a nonnegative integer,
  $$
S_{\mu}(z)=\exp\Big(-\int_0^{2\pi}\frac{e^{it}+z}{e^{it}-z}\,
d\mu(t)\Big) 
    $$
with positive singular measure $d\mu$  and
  $$
F(z)=\lambda\exp\Big(\frac{1}{2\pi}\int_0^{2\pi}\frac{e^{it}+z}{e^{it}-z}
\log\left| F^*\big(e^{it}\big)\right|\,dt\Big),\eqno (2)
 $$
where $\lambda\in\Bbb C$ with  $|\lambda|=1$ and  $\log\vert F^*\vert\in L^1(\Bbb T)$.\par  

Recall that a function $I$ of the form 
 $$
I(z)=B(z)S_{\mu}(z),\quad z\in \Bbb D,
  $$ 
is called an {\it inner function}. Furthermore, it is well known that
 $|I^*(e^{it})|=1$
for almost every  $e^{it}\in\Bbb T$   and hence,
$|f^*(e^{it})|=|F^*(e^{it})|$  for almost every  $e^{it}\in\Bbb T$.

I.I. Privalov \cite[p. 98]{p}  (alo see \cite[Theorem 5.3]{mp8}) proved 
that a function $f$ holomorphic on $\Bbb D$ belongs 
to the class   $N^p$ if and
only if $f=IF$, where $I$ is an inner function on $\Bbb D$ and $F$ is an outer
function given by (2) such that $\log^+\vert f^*\vert\in L^p(\Bbb T)$
(or equivalently, $\log^+\vert F^*\vert\in L^p(\Bbb T)$).

M. Stoll \cite[Theorem 4.2]{s} showed that the space $N^p$ (with the notation
 $(\log^+H)^\alpha$ in \cite{s}) equipped with the topology given by the metric 
$\lambda_p$ defined by            
   $$
\lambda_p(f,g)=\Bigg(\int_0^{2\pi}\big(\log(1+
\vert f^*(e^{i\theta})-g^*(e^{i\theta})\vert)\big)^p\,\frac{d\theta}
{2\pi}\Bigg)^{1/p},\quad f,g\in N^p,\eqno (3) 
   $$  
becomes an $F$-algebra.
Recall that the function $\lambda_1=\lambda$ defined on the Smirnov 
class $N^+$
by (3)  with $p=1$ induces the metric topology on $N^+$.
N. Yanagihara \cite{Y2} proved that 
under this topology, $N^+$ is an $F$-space.

Furthermore,  in connection with the spaces $N^p$
$(1<p<\infty)$,  Stoll \cite{s} 
(also see \cite{E1} and \cite[Section 3]{mp3})   also studied the spaces
$F^q$ $(0<q<\infty)$ (with the notation $F_{1/q}$ in \cite{s}), 
consisting of those functions $f$ holomorphic on 
$\Bbb D$ for which
   $$
 \lim_{r\to 1}(1-r)^{1/q}\log^+M_{\infty}(r,f)=0,
   $$ 
where
$$
 M_{\infty}(r,f)=\max_{\vert z\vert\le r}\vert f(z)\vert.
 $$
 Stoll  \cite[Theorem 3.2]{s} also proved that the space $F^q$ 
with the topology given by the  family of seminorms 
$\left\{|\Vert \cdot\Vert  |_{q,c}\right\}_{c>0}$ defined for $f\in F^q$ as 
  $$
|\Vert f\Vert|_{q,c}=\sum_{n=0}^{\infty}|\hat{f}(n)|e^{-cn^{1/(q+1)}}<\infty
  $$ 
for each $c>0$, where  $\hat{f}(n)$ is the $n$-th Taylor coefficient of $f$,
is a countably normed {\it Fr\'{e}chet algebra}. 
By a result of C.M. Eoff \cite[Theorem 4.2]{E1}, $F^p$ is the
{\it Fr\'{e}chet envelope} of $N^p$ and hence,
$F^p$ and $N^p$ have the same topological duals.

Following H.O. Kim (\cite{k1} and \cite{k2}), 
the class $M$ consists of all holomorphic functions $f$
on $\Bbb D$ for which
\begin{equation*}
\int_0^{2\pi}\log^+Mf(\theta)\frac{d\theta}{2\pi}<\infty,
\end{equation*}
where 
\begin{equation*}
Mf(\theta)=\sup_{0\leqslant r<1}\big\lvert f\big(re^{i\theta}\big)\big\rvert
\end{equation*}
is the {\it maximal radial function} of $f$.

The study on the class $M$ on the disk $\Bbb D$ has been exstensively
investigated by H.O. Kim in \cite{k1} and \cite{k2}, 
V.I. Gavrilov and V.S. Zaharyan \cite{gz} 
and M. Nawrocky \cite{n}.  
Kim \cite[Theorems 3.1 and 6.1]{k2} showed that the space $M$ with 
the topology given by the metric $\rho$ defined by
    $$
\rho(f,g)=\int_0^{2\pi}\log(1+M(f-g)(\theta))\,\frac{d\theta}{2\pi},
\quad f,g\in M,\eqno (4)
    $$    
becomes an $F$-algebra.
Furthermore,  Kim \cite[Theorems 5.2 and 5.3]{k2} gave
an incomplete characterization of multipliers of $M$ into 
$H^\infty$. 
Consequently, the topological dual of $M$
is not exactly determined in \cite{k2}, 
but as an application, it was proved in \cite[Theorem 5.4]{k2}
(also cf. \cite[Corollary 4]{n}) that  $M$ is not locally convex space.
Furthermore, the space $M$ is not locally bounded (\cite[Theorem 4.5]{k2}
and \cite[Corollary 5]{n}).

Nevertheless that as noticed above, the class $M$ 
is essentially smaller than the class $N^+$,
 M. Nawrocky \cite{n} showed 
 that the class $M$ and the Smirnov class $N^+$ have the same corresponding 
locally convex structure which was already established by 
N. Yanagihara for the Smirnov class in \cite{Y2} and  \cite{Y1}.
More precisely, it was proved in \cite[Theorems 1]{n}  
that the Fr\'{e}chet envelope of the  class $M$ can be 
identified with the space $F^+$
of holomorphic functions on the open unit disk $\Bbb D$ such that 
  $$
|\lVert f\rVert|_c:=\sum_{n=0}^{\infty}|\hat{f}(n)|e^{-c\sqrt{n}}<\infty
  $$ 
for each $c>0$, where  $\hat{f}(n)$ is the $n$-th Taylor coefficient of $f$.
Notice that $F^+$ coincides with  the space $F^1$ 
defined above.  
It was shown in \cite{Y1} (also see \cite{Y2}) that $F^+$ is actually 
the containing Fr\'{e}chet space for $N^+$ (also see \cite{St2}). 
Moreover, Nawrocky \cite[Theorem 1]{n} characterized the set of all 
continuous linear functionals on $M$  which by a result of Yanagihara \cite{Y2}
coincides with those on the Smirnov class $N^+$.  

Motivated by the mentioned investigations of the classes
$M$ and $N^{+}$, and the fact that the classes
$N^p$ $(1<p<\infty)$ are generalizations 
of the Smirnov class $N^+$, in \cite[Chapter 6]{me2} and \cite{me3}
 the first author of this 
paper  investigated  the classes  $M^p$ $(1<p<\infty)$
as generalizations of the class $M$. Accordingly, the {\it class} $M^p$
$(1<p<\infty)$ consists of all holomorphic functions $f$ on 
$\Bbb D$ for  which 
    $$
\int_0^{2\pi}\left(\log^+Mf(\theta)\right)^p\,\frac{d\theta}{2\pi}<\infty.
    $$
Obviously,     
   $$
\bigcup_{p>1}M^p\subset M.
   $$
By analogy with the topology defined on the  space $M$
(\cite{k1} and \cite{k2}), the space $M^p$ can be  equipped with the topology 
induced by the metric $\rho_p$ defined as            
 $$
\rho_p(f,g)=\Big(\int_0^{2\pi}\log^p(1+M(f-g)(\theta))\,
\frac{d\theta}{2\pi}\Big)^{1/p},
 $$
with $f,g\in M^p$.

After Privalov, the study of the spaces $N^p$ $(1<p<\infty)$ was continued in 1977
by M. Stoll  \cite{s} (with the notation
 $(\log^+H)^\alpha$ instead of $N^p$ in \cite{s}).  Further, the 
linear topological and 
functional properties of these spaces  were extensively investigated 
 by C.M. Eoff in \cite{E1} and  \cite{E2}, N. Mochizuki \cite{Mo},
Y. Iida and N. Mochizuki \cite{im}, Y. Matsugu \cite{ma},
J.S. Choa \cite{ch}, J.S. Choa and H.O. Kim \cite{CK},
 A.K. Sharma and S.-I. Ueki  \cite{su}
 and in works \cite{me4}-\cite{msu}  of  authors of this paper; typically, 
the notation varied and Privalov was mentioned in \cite{ma},
\cite{me5}-\cite{me7},  \cite{mp3}-\cite{mp7}, 
\cite{ms}, \cite{msu} and \cite{su}. 
In particular, it was proved in \cite[Corollary]{me5} that  $N^p$ 
is not locally convex space  and  in  \cite[Theorem 1.1]{mp4} that 
$N^p$ is not locally bounded space. 
We refer the recent monograph \cite[Chapters 2, 3 and 9]{gse} by  
V.I. Gavrilov, A.V. Subbotin and D.A. Efimov  for a good reference on the 
spaces $N^p$ $(1<p<\infty)$. 

Let us recall that in our recent  overview paper \cite{mp7} it was given  a 
survey of some known results on different topologies  on the Privalov 
classes $N^p$ $(1<p<\infty)$ and their Fr\'{e}chet envelopes 
$F^p$ $(1<p<\infty)$ on the open unit
disk. Here we give a survey on related extended results involving 
some other metrics and the induced  topologies on the classes 
$N^p$. 

The remainder of this overview paper is organized in three sections.
For any fixed $p>1$, in Section 2  we present 
some results concerning the topological and functional structures 
on the classes $M^p$ $(1<p<\infty)$. Section 3 is devoted to the consideration 
of the Privalov class $N^p$ as a closed subspace of some Orlicz space.
In this setting $N^p$ with the associated modular in the sense 
of Musielak and Orlicz  
becomes the Hardy-Orlicz class whose topology coincides with 
both metric topologies $\lambda_p$ and $\rho_p$. 
Concluding remarks are presented in the last section.

\vspace{6mm}

\noindent {\bf  2 THE $\rho_p$-METRIC TOPOLOGY ON PRIVALOV SPACE $N^p$}

\vspace{3mm}

Here we focus our attention to certain results from  \cite[Chapter 6]{me2} 
and \cite{me3} concerning the  classes $M^p$ $(1<p<+\infty)$.
In \cite{me3} it is proved   the following basic result. 

\vspace{3mm}
\noindent{\bf Theorem 1} (\cite[Theorem 2]{me3}).
{\it The function $\rho_p$ defined on $M^p$ as
    $$
\rho_p(f,g)=\Big(\int_0^{2\pi}\log^p(1+M(f-g)(\theta))\,
\frac{d\theta}{2\pi}\Big)^{1/p},\quad f,g\in M^p, \eqno(5)
     $$
is a translation invariant metric on $M^p$. Further, the space  $M^p$ 
is a complete metric space with respect to the metric $\rho_p$.}
  \vspace{3mm}

\noindent{\bf Remark 1.} Notice that the expression (5) with  
$p=1$ defines the metric  $\rho_1=\rho$ on the class 
$M$ (given by (4)) introduced   by H.O. Kim in \cite{k1} and \cite{k2}.  
As noticed above, it was proved in \cite{k2} that the metric  $\rho$ 
induces the topology on $M$ under which $M$ is also an $F$-algebra. 
  \vspace{3mm}

Moreover, the following two statements  are  also proved in 
\cite[Chapter 6]{me2}. 

\vspace{3mm}

\noindent{\bf Theorem 2}  (\cite[Theorem 11]{me3}).  
{\it  $M^p=N^p$ for each $p>1$; that is, the spaces $M^p$ and $N^p$ 
coincide.}

\vspace{3mm}

\noindent{\bf Theorem 3} (\cite[Theorem 15]{me3}).  
{\it $M^p$ with the topology given by the metric 
$\rho_p$ defined by {\rm (5)}  becomes an $F$-space.}

\vspace{3mm}

Using Theorem 3 and  the open mapping theorem 
(see, e.g., \cite[Corollary  2.12 (b)]{r}), the following result was also 
proved in \cite{me3}.

\vspace{3mm}

\noindent{\bf Theorem 4} (\cite[Theorem 16]{me3}). 
{\it For each $p>1$ the classes  $M^p$ and $N^p$ coincide, 
and the metric spaces $(M^p,\rho_p)$ and $(N^p,\lambda_p)$ have the same 
topological structure, where the metrics $\rho_p$ and 
$\lambda_p$ are given on $M^p$ and $N^p$ by $(5)$ and $(3)$, respectively.}

\vspace{3mm}

As an immediate consequence of Theorem 4 and \cite[Lemma 8]{me3},
we  obtain the following assertion. 

\vspace{3mm}

\noindent{\bf Proposition 1.} 
{\it The convergence
with respect to the metric $\rho_p$ given by $(5)$ on the space $M^p$  
is stronger than the metric of uniform convergence
on compact subsets of the disk $\Bbb D$.}

\vspace{3mm}

\noindent{\bf Remark 2.} For an outer function $h$ let 
$H^2(|h^*|^2)$ denote the closure of the (analytic) polynomials
in the space $L^2(|h^*|^2d\,\theta)$.  By using the famous   
Beurling's theorem for the Hardy space 
$H^2$ (\cite{b}; also see 
\cite[Ch. 7, p. 99]{ho}), it was proved in  \cite{E2}
(also see \cite[Section 1]{mp2}) that the class  $N^p$ can be 
represented as a union of certain weighted Hardy classes.     
Using this representation, the following 
two topologies are  defined on $N^p$ in \cite{mp2}:
the usual locally convex inductive limit topology, which we shall call the
 {\it Helson topology} and denote by ${\mathcal H}_p$, in which a neighborhood base for $0$ 
 is given by those balanced convex sets whose intersection with each  
$H^2\left(\vert h^*\vert^2\right)$ is a neighborhood of zero in $H^2\left(
\vert h^*\vert^2\right)$, and a not locally convex topology, denoted
by $I_p$, in which a neighborhood base for zero is given by all sets
whose intersection with each space $H^2\left(\vert h^*\vert^2\right)$ is a
neighborhood of zero. It was proved in 
 \cite[Theorem E]{mp2} (cf. \cite[Chapter 3]{me2}) that the  topology 
${\mathcal H}_p$  
coincides with the metric topology induced on $N^p$ by the 
Stoll's metric topology  $\lambda_p$ given by  (3). 
Moreover,  it was proved in \cite{E2} 
that the topology $I_p$ coincides 
with the metric topology $\lambda_p$ and by Theorem 4, $I_p$ also coincides 
 with the metric topology $\rho_p$, which are not locally convex. 
Hence, $I_p$ is strictly stronger than ${\mathcal H}_p$. 
The analogous results for the space
$N^+$ are proved by J.E. McCarthy in  \cite{Mc2}.

\vspace{6mm}

\noindent{\bf 3 THE SPACE  $N^p$ AS THE HARDY-ORLICZ CLASS}
  \vspace{3mm}

In this section we we give a short survey about 
 Privalov  classes $N^p$ $(1<p<+\infty)$  
as the  Hardy-Orlicz classes.
Related results are mainly obtained in \cite{mpl}. 
Let $(\Omega, \Sigma, \mu)$ be a measure space, i.e., $\Omega$ is a
nonempty set, $\Sigma$ is a $\sigma$-algebra of subsets of $\Omega$ and
$\mu$ is a nonnegative finite complete measure not vanishing identically.
Denote by $ L^p(\mu)= L^p$ $(0<p\le\infty)$ the familiar 
Lebesgue spaces on
$\Omega$. For each real number $p>0$ in \cite{mpl}  it was  considered  
 the class 
$L^+_p(\mu)={L}^+_p$ of all (equivalence classes of) 
$\Sigma$-measurable complex-valued functions
$f$ defined on $\Omega$ such that the function $\log^+|f|$
belongs to the space $L^p$,  i.e., 
 $$
\int_{\Omega}\left(\log^+|f(x)|\right)^p\,
d\mu<+\infty,
 $$
where $\log^+|a|=\max (\log |a|,0)$.
 Clearly, $L^+_q\subset L^+_p$ for $q>p$ 
and  
$\bigcup_{p>0}L^p\subset\bigcap_{p>0}L^+_p$ \cite[Section 2]{mpl}. 
For each $p>0$ the  space ${L}^+_p$ is an algebra 
with respect to the pointwise addition and multiplication.
 For each $p>0$ we define the metric $d_p$ on $L^+_p$ by
\begin{eqnarray*}
\qquad \qquad d_p(f,g)&=&\inf_{t>0}\left[t+\mu\left(\{x\in\Omega
:\,|f(x)-g(x)|\ge t\}\right)\right]\\
\qquad \qquad&&+\int_{\Omega}\left|\left(\log^+|f(x)|
\right)^p-\left(\log^+|g(x)|\right)^p\right|\,d\mu.\qquad
\qquad\qquad\quad(6)
\end{eqnarray*} 
Recall that the space $L^+_1$ was introduced by T. Gamelin and 
 G. Lumer in \cite[p. 122]{gl} (also see  
 \cite[p. 122]{g}, where $L^+_1$ is denoted as  $L(\mu)$). Note that   the  metric $d_p$ given by (6) 
with $p=1$ coincides with the
Gamelin-Lumer's  metric $d$ defined on $L^+_1$. It was proved in 
\cite[Theorem 1.3, p. 122]{gl}  (also see \cite[Theorem 2.3, p. 122]{g})
 that the space
$L^+_1$ with the topology given by the metric $d_1$ becomes a topological 
algebra.
The following result is a generalization of the corresponding result
for the case $p=1$ given in \cite[p. 122]{gl} (also see \cite[p. 122]{g}).

\vspace{3mm}
\noindent{\bf Theorem 5} (\cite[Theorem 2.1]{mpl}).
{\it \vspace{3mm}
The space $L^+_p$ with the metric $d_p$ given by $(6)$ is a  topological 
algebra, i.e., a topological vector space with a complete metric in which
multiplication is continuous.}
\vspace{3mm}

By the inequality 
 $$
\big(\log(1+|z|)\big)^p\le 2^{\max (p-1,0)}\big((\log 2)^p+(\log^+|z|)^p\big),
\quad z\in\Bbb C,
 $$ 
it follows that a function $f$ belongs to the space ${L}^+_p$
if and only if 
 $$
\|f\|_p:=\Bigg(\int_{\Omega}\left(\log(1+|f(x)|)\right)^p\,d\mu\Bigg)^{1/p}
<\infty.\eqno(7)
 $$
Furthermore \cite[Section 2]{mpl}, the function $\sigma_p$ defined  as
 $$
\sigma_p(f,g)=\big(\|f-g\|_p\big)^{\min (1,p)},\quad f,g\in L^+_p, 0<p\le 1,
\eqno(8)
 $$
is a translation invariant metric on $L^+_p$ for all $p>0$.
 Notice that in the case of Privalov space $N^p$ $(1<p<\infty)$,
the metric $\sigma_p$ given by (8) coincides with Stoll's metric 
$\lambda_p$ defined by (3).

Recall that two metrics (or norms) defined on the same space will be called 
 equivalent if they induce the same topology on this space.

\vspace{3mm}
\noindent{\bf Theorem 6} (\cite[Theorem 2.3]{mpl}).
{\it The metric $d_p$ given by $(6)$  defines the topology for $L^+_p$ which is
equivalent to the topology defined by the metric $\sigma_p$ given by 
$(8)$}. 
\vspace{3mm}

\noindent{\bf Remark 4.} 
It was pointed out in \cite[Remark, Section 2]{mpl} that using the same argument 
applied in the proof of Theorem 2.3 of \cite{mpl}, 
it is easy to show that the metrics $\sigma_p$ and $d_p$
are equivalent with the metric $\delta_p$ given on $L^+_p$ by
\begin{eqnarray*}
\qquad \delta_p(f,g)&=&\inf_{t>0}\left[t+\mu\left(\{x\in\Omega
:\,|f(x)-g(x)|\ge t\}\right)\right]\\
&&+\Bigg(\int_{\Omega}\left|\log^+|f(x)|
-\log^+|g(x)|\right|^p\,d\mu\Bigg)^{1/\max(p,1)},
f,g\in L^+_p.  \quad\,\,\, (9)
\end{eqnarray*}
\vspace{3mm}

\noindent{\bf Remark 5.} 
In  \cite[Remark  5, p. 460]{Y2} M. Hasumi pointed out
that the Yanagihara's metric $\lambda=\lambda_1$ on the Smirnov class
 (given by (3) with $p=1$)  defines the topology on 
the space $L^+_1=L(\mu)$ which is equivalent to the metric  topology  
$d_1=d$ (given by (6) with $p=1$).
 \vspace{3mm}

As a consquence of Theorems 5 and 6, it can be obtained the following 
result.

\vspace{3mm}
\noindent{\bf Theorem 7} (\cite[Corollary 2.4]{mpl}).
{\it For each $p>0$ the space $L^+_p$ with the topology given by the
metric $\sigma_p$ is an $F$-algebra, i.e., a topological algebra with a 
complete translation invariant metric $\sigma_p$.}  
\vspace{3mm}

\noindent{\bf Remark 6.}
In view of Theorem 7, note that
$L^+_p$ may be considered as the generalized Orlicz space  $L^w_p$
with the constant function $w(t)\equiv 1$ on $[0,2\pi)$ 
defined in \cite[Section 6]{mpl}.

\vspace{3mm}

The real-valued function $\psi :[0,\infty)\mapsto [0,\infty)$
defined as $\psi(t)=\big(\log(1+t)\big)^p$,
is continuous and nondecreasing in $[0,\infty)$, equals zero only 
at 0, and hence it  is a $\varphi$-{\it function} (see, e.g.,  
\cite[p. 4, Examples 1.9]{mu}). Moreover, $\psi$ 
is a log-{\it convex} function since it can be represented in the form 
$\psi (x)=\Psi (\log x)$ for $x>0$, where $\Psi (u):=\max (u^p,0)$
($u\in [0,\infty)$) is a {\it convex} function on the whole real axis,
satisfying the condition $\lim_{u\to +\infty}\frac{\Psi(u)}{u}=+\infty$.
Notice that convex $\varphi$-functions are a particular case of log-convex 
$\varphi$-functions. 

Further,
observe that \cite[Section 4]{mpl} the space
 $L^+_p(dt/(2\pi)=L_p^{+}$  $(p>0)$, consisting of all 
complex-valued functions $f$, defined and measurable on $[0,2\pi)$, for which
 $$
\|f\|_p:=\left(\int_{0}^{2\pi}\left(\log(1+|f(t)|)\right)^p\,
\frac{dt}{2\pi}\right)^{1/p}<+\infty\eqno(10)
 $$
is the {\it Orlicz class} (see \cite[p. 5]{mu}; 
cf. \cite[Section 4]{mpl}), whose generalization 
was given in  \cite[Section 6]{mpl}.
It follows by the dominated convergence theorem
that the class $L_p^{+}$ coincides with the
associated  {\it Orlicz space} (see \cite[Definition 1.4, p. 2]{mu}) 
consisting 
of those functions $f\in L_p^{+}$ such that
\[
\int_{0}^{2\pi}\left(\log(1+c|f(t)|)\right)^p\,
\frac{dt}{2\pi}\to 0\quad {\mbox as }\quad c\to 0+.
\]
Since  $\sigma_p(f,g)=\big(\|f-g\|_p\big)^{\min (p,1)} (f,g\in L_p^{+}$)
is an invariant metric on $L_p^{+}$, 
the function $\|\cdot\|_p$ given by (10)
is a {\it modular} in the sense of Definition 1.1 in \cite[p. 1]{mu},
where $\sigma_p$ is the metric defined by (8).  
For any function  $f\in L_p^{+}$, by the monotone convergence theorem, 
it follows that $\lim_{c\to 0}\Vert cf\Vert =0$
and thus $(L_p^{+},\sigma_p)$ is a {\it modular
space} in the sense of Definition 1.4 in \cite[p. 2]{mu}.
 In other words,
the function $\|\cdot\|_p$  is an $F$-norm. It is  known (see  
\cite[Theorem 1.5, p. 2 and Theorem 7.7, p. 35]{mu}) 
that the functional $|\cdot|_p$ defined as
 \[
|f|_p=\inf\Big\{\varepsilon>0:\,\int_{0}^{2\pi}\left(\log\left(1+\frac{|f(t)|}
{\varepsilon}\right)\right)^p\,\frac{dt}{2\pi}\le\varepsilon
\Big\},\quad f\in L_p^{+},\eqno(11)
  \]          
is a complete $F$-{\it norm} on $L_p^{+}$. Furthermore (see \cite[p. 54]{l}), 
 if we denote by $L_p^{0}$ the class of all 
functions $f$ such that $\alpha f\in L_p^{+}$ for  
every $\alpha>0$, then $L_p^{0}$   is the 
closure of the space of all continuous functions on $[0,2\pi)$
in the space $(L_p^{+},|\cdot|_p)$.  

We also give the following two results.

\vspace{3mm}

\noindent{\bf Theorem 8} (\cite[Theorem 4.1]{mpl}).
{\it \vspace{3mm}
The $F$-norms $\|\cdot\|_p$ and $|\cdot|_p$
$($given by $(10)$ and $(11)$, respectively$),$  induce the same 
topology on the space $L_p^{+}$. In other words, the norm and modular convergences
are equivalent.}

\vspace{3mm}

\noindent{\bf Proposition 2} (\cite[Corollary 4.2]{mpl}).
{\it \vspace{3mm} There does not exist a nontrivial continuous linear functional on the space
$(L_p^{+},\|\cdot\|_p)$.}
\vspace{3mm}

Note  that \cite[Section 5]{mpl} the algebra $N^p$ may be considered as the 
{\it Hardy-Orlicz
space} with the Orlicz function $\psi :[0,\infty)\mapsto [0,\infty)$
defined as $\psi(t)=\big(\log(1+t)\big)^p$. These spaces were firstly 
studied in 1971 by  R. Le\'{s}niewicz \cite{l2}. 
For more information on the
Hardy-Orlicz spaces, see \cite[Ch. IV, Sec. 20]{mu}. Identifying a function
$f\in N$ with its {\it boundary function} $f^*$, by \cite[3.4, p. 57]{l}, the 
space $N^p$ is identical with the closure of the space of all functions
holomorphic on the open unit disk  $\Bbb D$ and continuous on 
 $\overline{\Bbb D}:|z|\le 1$ 
in the space
$(L_p^{+}(dt/2\pi)\cap N,|\cdot|_p)$, where $dt/2\pi$ is the 
usual normalized Lebesgue measure on the unit circle $\Bbb T$.   
Using this fact, Theorem 4 and  Theorem 6, the main  surveyed results of 
this paper can be summarized as follows.

 \vspace{5mm}
\noindent{\bf Theorem 9} (\cite[Theorem ]{mpl}).
{\it \vspace{3mm} For each $p>1$ Privalov class 
is the Hardy-Orlicz space with the Orlicz 
function  $\psi(t)=\big(\log(1+t)\big)^p$ $($$t\in [0,2\pi)$$)$.   
Moreover, the metrics $\lambda_p$, $\rho_p$, $d_p$, $\delta_p$
and the functional $|\cdot |_p$ $($defined respectively 
by} (3), (5), (6), (9) {\it and} (11)) {\it induce the same 
topology on $N^p$ under which $N^p$ becomes an $F$-algebra.} 
  \vspace{5mm}

\vspace{6mm}

\noindent {\bf  4 CONSLUSION}

\vspace{2mm}

This paper continues an overview of topologies on the Privalov spaces
$N^p$  ($1<p<+\infty$)  induced by different metrics.
Notice that the class $N^p$
equipped with the topology given by the metric $\lambda_p$
introduced by M. Stoll becomes an $F$-algebra.
The same statement is also true for the class $M^p$ 
with respect to the $\rho_p$-metric topology.
These  facts are used in \cite{me2} to prove  that for each 
$p>1$  the classes  $M^p$ and $N^p$ coincide 
and the metric spaces $(M^p,\rho_p)$ and $(N^p,\lambda_p)$ 
have the same topological structure
In Section 3 we  give a short survey 
about  Privalov  spaces $N^p$ $(1<p<+\infty)$  
whose topology is induced by the generalized Gamelin-Lumer's 
metric $d_p$  defined on the
space $L^+_p(dt/(2\pi)$. Notice that  the space $L^+_p(dt/(2\pi)$ 
coincides with the   Orlicz class associated to the 
log-convex $\varphi$-function $\psi(t)=\big(\log(1+t)\big)^p$
($t\in [0,+\infty $). Accordingly, it follows that
for each $p>1$ Privalov space
is the Hardy-Orlicz space with the Orlicz 
function  $\psi(t)=\big(\log(1+t)\big)^p$ ($t\in [0,2\pi)$).    
Moreover, the metrics $\lambda_p$, $\rho_p$ and $d_p$
induce the same 
topology on $N^p$ under which $N^p$ becomes an $F$-algebra.
We believe that presented results would be useful for future  
research on related topics, as well as for some applications in 
Functional and Complex Analysis.  
      
 \vspace{2mm}

{\noindent{\bf Acknowledgements.}} This work was supported by the MEPhI Academic 
Excellence Project (agreement with the Ministry of Education and Science of 
the Russian Federation of August 27, 2013, project no. 02.a03.21.0005).

\end{document}